\documentclass[letterpaper,11pt]{amsart}

\usepackage{color}
\usepackage[T1]{fontenc}
\usepackage{amsthm}
\usepackage{amsmath}
\usepackage{amssymb}
\usepackage{amssymb, amsfonts, amsbsy, latexsym, epsfig, color, verbatim}
\usepackage{enumerate}

\newtheorem{thm}{Theorem}[section]

\newtheorem{cor}[thm]{Corollary}

\newtheorem{problem}[thm]{Problem}
\newtheorem{lem}[thm]{Lemma}

\newtheorem{prop}[thm]{Proposition}
\newtheorem{dfn}[thm]{Definition}

\newtheorem{notation}[thm]{Notation}

\theoremstyle{remark}

\newtheorem{example}[thm]{Example}

\newtheorem{rmk}[thm]{Remark}

\newcommand{\ip}[2]{\ensuremath{\left\langle #1, #2 \right\rangle}}
\newcommand{\absip}[2]{\ensuremath{\left| \left\langle #1, #2 \right\rangle \right|}}

\DeclareMathOperator\diag{diag}

\title{Integer Frames}
\author[P.G. Casazza, R.G. Lynch, J. Tremain, and L.M. Woodland
 ]{Peter G. Casazza, Richard G. Lynch, Janet C. Tremain, and Lindsey M. Woodland }
\address{Department of Mathematics, University
of Missouri, Columbia, MO 65211-4100}

\thanks{The authors were supported by NSF 1307685; NSF ATD 1042701; NSF ATD 00040683; AFOSR DGE51: FA9550-11-1-0245}

\email{casazzap@missouri.edu, rglz82@mail.missouri.edu,\newline tremainjc@missouri.edu, lmwvh4@mail.missouri.edu}

\begin{document}
\begin{abstract}
Finite frame theory has become a powerful tool for many applications of mathematics.
  In this paper we introduce a new area of research in frame theory:
Integer frames.  These are frames having all integer coordinates with
respect to a fixed orthonormal basis for a Hilbert space.  Integer
frames have potential to mitigate quantization errors
and transmission losses as well as speeding up computation times.  This paper gives
the first systematic study of this important class of finite Hilbert space frames.
\end{abstract}

\maketitle

\section{Introduction}

{\it Integer frames}, which are frames whose vectors have all integer coordinates
with respect to a fixed orthonormal basis for a Hilbert space,
have the potential to mitigate quantization errors and transmission loses as well as
speed up computation time.  In this paper we give the first systematic study of
this class of frames.
The focus of the present paper is the construction of such frames.  The main
goal is to give construction methods for integer frames
with the added properties of equal norm, tight, and/or full spark.  However, dropping either one of the assumptions that the frame be equal norm or tight is also considered.

We note that to construct integer frames, it suffices to construct frames with
rational coordinates because we can then multiply the frame by the greatest common denominator of the rationals in order to get an integer frame.

The paper is arranged as follows.
Section 2 gives the necessary background material from finite frame theory used throughout the paper. Section 3 covers three propositions that give a method of constructing larger frames from those with fewer vectors or those in lower dimensions. An application concerning {\it Hadamard} matrices is also discussed. Section 4 deals with equal norm, tight, integer frames in two and three dimensions, in which two dimensions is answered completely and only partial results are given for three dimensions. In Section 5, the special case of frames having one more element than the dimension is examined. It is shown that the existence of frames having $M+1$ vectors in $M$ dimensions is directly related to the existence of $M$-simplexes having integer coordinates in $M$ dimensions. Section 6 is dedicated to equal norm, tight frames in general dimension. Finally, it is shown in Section 7 that when dropping either the equal norm or tight assumptions, any number of vectors in any dimension can be obtained. The same is shown for equal norm frames that are {\it nearly tight}.

\section{Frame Theory}

A brief introduction to frame theory is given in this section, which contains the necessary background for this paper. For a thorough approach to the basics of frame theory see \cite{petesbook, ole_book}.

\begin{dfn}
A family of vectors $\{f_i\}_{i=1}^N$ in an $M$-dimensional Hilbert space, $\mathcal{H}_M$, is a \emph{frame} if there are constants $0 < A \leq B < \infty$ so that for all $f \in \mathcal{H}_M$,
\[
A \|f \|^2 \leq \sum_{i=1}^N \absip{f}{f_i}^2 \leq B \|f\|^2,
\]
where $A$ and $B$ are the \emph{lower frame bound} and \emph{upper frame bound}, respectively. If $A = B$, this is a \emph{tight frame} and if $A = B = 1$, it is a \emph{Parseval frame}. If there is a constant $c$ so that $\|f_i\| = c$ for all $i = 1,\ldots,N$, it is an \emph{equal norm frame} and if $c = 1$, then it is a \emph{unit norm frame}. If there is a constant $d$ so that $|\ip{f_i}{f_j}| = d$ for all $i \neq j$, then it is an \emph{equiangular frame}. Finally, the values $\{ \ip{f}{f_i} \}_{i = 1}^N$ are called the \emph{frame coefficients} of the vector $f\in \mathcal{H}_M$ with respect to the frame $\{f_i\}_{i=1}^N$.
\end{dfn}

If $\{f_i\}_{i=1}^N$ is a frame for $\mathcal{H}_M$, the \emph{analysis operator} of the frame is the operator $T : \mathcal{H}_M \to \ell_2(N)$ given by
\[
T(f) = \{ \ip{f}{f_i} \}_{i=1}^N
\]
and the associated \emph{synthesis operator} is the adjoint operator $T^*$ and satisfies
\[
T^{*}\left( \{a_i\}_{i=1}^N \right) = \sum_{i=1}^N a_i f_i.
\]
The \emph{frame operator} is the positive, self-adjoint, invertible operator $S = T^* T$ on $\mathcal{H}_M$ and satisfies
\[
S(f) = T^*T(f) = \sum_{i=1}^N \ip{f}{f_i} f_i, \hspace{.2in} f \in \mathcal{H}_M.
\]
If $S$ has eigenvalues $\{\lambda_j\}_{j=1}^M$, then
\[
\sum_{i = 1}^N \|f_i\|^2 = \sum_{j = 1}^M \lambda_j
\]
and the largest and smallest eigenvalues of the frame operator $S$ coincide with the optimal upper frame bound and the optimal lower frame bound, respectively.

In finite dimensions a frame is simply a spanning set, and as such the zero vector could potentially be one or more elements of a frame. However, since we are concerned with using integer frames to mitigate quantization errors and to speed up computation time, then the zero vector is not useful in our application since it provides no new information. Because of this, in the present paper, we will assume that no frames contain the zero vector.

It is important to note that for any frame $\{f_i\}_{i=1}^N$ in $\mathcal{H}_M$ with analysis operator $T$, the matrix representation of the synthesis operator $T^*$ with respect to some orthonormal basis $\{e_i\}_{i=1}^M$ of $\mathcal{H}_M$ is given by the following $M \times N$ matrix
$$\left[\begin{array}{cccc}
|&|&\dots&|\\
f_1&f_2&\dots&f_N\\
|&|&\dots&|
\end{array}\right]$$
where the columns of $T^*$ represent the coefficients of the frame vectors with respect to $\{e_i\}_{i=1}^M$. Due to this relationship between a frame and its matrix representation, we will not distinguish between a frame and its matrix and instead use the term \emph{frame} interchangeably.

Also, since the columns of the synthesis matrix represent the frame vectors, then the square sum of each column represents the square norm of the frame vectors. Hence a frame is equal norm if all of the columns square sum to the same constant.

This paper will be concerned with the construction of integer frames and will be approached by finding a matrix representation for the synthesis operator, which has all integer entries. Note that
for any rank $M$,
an $M \times N$ matrix with all integer entries represents an $N$-element integer frame in $\mathcal{H}_M$, where the frame vectors are the columns of this matrix with respect to an orthonormal basis for $\mathcal{H}_M$.

However,
for an arbitrary rank $M$, an $M \times N$ integer matrix, in general, does not have enough ``nice" properties to prove to be useful in applications. In applications, integer frames have the potential to speed up computation time and because of this, the frame operator and the eigenvalues of the frame operator should be readily available. The following theorem addresses this issue and defines the added properties needed when constructing ``application ready" integer frames.

\begin{thm} \label{matrix} \cite{petesbook}
Let $T : \mathcal{H}_M \to \ell_2(N)$ be a linear operator, let $\{e_j\}_{j=1}^M$ be an orthonormal basis for $\mathcal{H}_M$, and let $\{\lambda_j\}_{j=1}^M$ be a sequence of positive numbers. Let $A$ denote the $M \times N$ matrix representation of $T^{*}$ with respect to $\{e_j\}_{j=1}^M$ and the standard basis $\{\hat{e}_i \}_{i=1}^N$ of $\ell_2(N)$. Then the following conditions are equivalent.
\begin{enumerate}
\item $\{T^* \hat{e}_i \}_{i=1}^N$ forms a frame for $\mathcal{H}_M$ whose frame operator has eigenvectors $\{e_j\}_{j=1}^M$ and associated eigenvalues $\{\lambda_j\}_{j=1}^M$.
\item The rows of $A$ are orthogonal and the $j$-th row square sums to $\lambda_j$.
\item The columns of $A$ form a frame for $\ell_2(M)$ and $$A A^* = \diag(\lambda_1,\ldots,\lambda_M).$$
\end{enumerate}
\end{thm}

As a result of Theorem \ref{matrix}, it is clear that if we impose the synthesis matrix of a frame to be represented against the eigenbasis of its frame operator, $S$, then the synthesis matrix will have orthogonal rows and the square sum of the rows will be the eigenvalues of the frame operator $S$. Moreover, a frame is tight if all eigenvalues of $S$ are equal. Hence, in this situation, a frame is tight if the square sum of all rows are equal. Thus, constructing an integer frame against the eigenbasis of its frame operator will ensure that our frame has ``nice" properties and will help to speed up computation time in application. Because of this, all integer frames (unless stated otherwise) in the present paper will be represented against the eigenbasis of their frame operator and hence this requires orthogonality between the rows.

We will adopt the following notation:

\begin{notation}
We will write ENTIF for an equal norm, tight integer frame.
\end{notation}

The next result is basic; but since we use it extensively throughout
the paper, we record it formally here.

\begin{prop}
If $A=(a_{ij})_{i=1,j=1}^{  M\ ,\  N}$ is a frame matrix and $I\subset \{1,2,\ldots,M\}$ then
$B= (a_{ij})_{i\in I,j=1}^{\ \ \ \ \ N}$ is also a frame matrix.
\end{prop}

Finally, the notion of spark is introduced, which is the measure of how resilient a frame is against erasures, so {\it full spark} is a desired property of a frame.

\begin{dfn}
The  \emph{spark} of a frame $\{f_i\}_{i=1}^N$ in $\mathcal{H}_M$ is the cardinality of the smallest linearly dependent subset of the frame. The frame is called \emph{full spark} if every $M$ element subset of the frame is linearly independent.
\end{dfn}

In general, it is very difficult to check the spark of a frame. Moreover, it is shown in \cite{ajm_spark} that determining if a matrix is full spark is NP-hard.

\section{Combining Frames}\label{sec_oldnew}

In this section, we will see how to combine existing frames to obtain frames with more
vectors.
 These results will be used throughout the paper.  The next proposition is also
clear; but we record it for future reference.

\begin{prop}\label{adjoinprop1}
Let $A$ be an $M \times N_1$ matrix and $B$ be an $M \times N_2$ matrix and suppose $A$ and $B$ both represent frames in $\mathcal{H}_M$ with $N_1$ and $N_2$ elements, respectively. Then the $M \times (N_1 + N_2)$ block matrix $[A,B]$ represents a frame with $N_1 + N_2$ elements in $\mathcal{H}_M$. Furthermore, if $A$ and $B$ are both tight frames then $[A, B]$ is also a tight frame. Lastly, if $A$ and $B$ are both of the same equal norm, then $[A,B]$ is also equal norm.
\end{prop}

It is easy to see via induction that the preceding proposition also holds for any number of frames over the same Hilbert space. One can also adjoin the matrices diagonally which requires only that $A$ and $B$ be frames.  This result is also clear
so we omit its proof.

\begin{prop}\label{adjoinprop2}
Suppose $A$ and $B$ are $M_1 \times N_1$ and $M_2 \times N_2$ matrices which represent frames in $\mathcal{H}_{M_1}$ and $\mathcal{H}_{M_2}$, respectively. Then the $(M_1 + M_2) \times (N_1 + N_2)$ block diagonal matrix
\[
C = \left[\begin{array}{cc} A & \mathbf{0} \\ \mathbf{0} & B \end{array}\right]
\]
represents an $N_1+N_2$ element frame in $\mathcal{H}_{M_1 + M_2}$. For $C$ to be tight, $A$ and $B$ need to have the same tightness factor and for $C$ to be equal norm, both $A$ and $B$ need to be equal norm with the same factor.
\end{prop}

The last proposition of this section gives a method for constructing a new frame having twice the dimension and twice the number of elements.

\begin{prop}\label{adjoinprop3}
If $A$ is an $M \times N$ matrix representing a frame in $\mathcal{H}_M$ and $c$ is a nonzero scalar, then the $2M \times 2N$ matrix
\[
B = \left[\begin{array}{rr} cA & cA \\ cA & -cA \end{array}\right]
\]
represents a frame in $\mathcal{H}_{2M}$. The frame $B$ is tight if $A$ is tight and $B$ is equal norm if $A$ is equal norm.
\end{prop}

To demonstrate the usefulness of Proposition \ref{adjoinprop3}, we consider building ENTIFs out of \emph{Hadamard matrices}.

\begin{dfn}
An $N \times N$ matrix $A$, having only $\pm 1$ as its entries and satisfying $A^T A = N\cdot I_{N\times N}$ is called a \emph{Hadamard matrix}.
\end{dfn}

We are interested in Hadamard matrices because if an $N \times N$ Hadamard matrix, $A$, exists then the $M \times N$ matrix formed by the first $M$ rows of $A$ is an $N$-element ENTIF in $M$ dimensions. Also note that a Hadamard matrix itself represents an ENTIF and so Proposition \ref{adjoinprop3}, with $c = 1$, implies that if an $N\times N$ Hadamard matrix exists, then there is also a Hadamard matrix of size $2^k N \times 2^k N$ for all $k \geq 0$. Thus a frame with $2^k N$ elements can also be formed in $M$ dimensions. This is summarized in the following theorem.

\begin{thm}\label{bighadthm}
Suppose an $N \times N$ Hadamard matrix exists for some $N \in \mathbb{N}$. Then an ENTIF with $2^k N$ elements in $M$ dimensions exists for all $k \geq 0$ and $M \leq 2^k N$.
\end{thm}

The preceding theorem is a generalization of a now standard construction of Sylvester, who showed that $2^K \times 2^K$ Hadamard matrices exist for all nonnegative integers $K$. Namely, let $H_0$ be the $1 \times 1$ matrix
\[
H_0 = \left[ 1 \right]
\]
and iterate to obtain the $2^K \times 2^K$ matrix
\[
H_K = \left[\begin{array}{rr} H_{K-1} & H_{K-1} \\ H_{K-1} & -H_{K-1} \end{array}\right]
\]
for any positive integer $K$. Now forming a new matrix by choosing the first $M \leq 2^K$ rows of $H_K$ yields an ENTIF with the square norms of the columns (frame vectors) equal to $M$. It is worth noting that the ENTIF obtained in this way may not be full spark since, for instance, keeping only the first half of the rows of $H_K$ to form a  frame gives two copies of an orthonormal basis. In general, it is not known which subsets of the rows of a Hadamard matrix give a full spark frame.

It is a well-known result that $N \times N$ Hadamard matrices can only exist when $N = 1,2,4K$, where $K \geq 1$. However, the existence of a Hadamard matrix of size $4K$ is not yet known for all values of $K$ and the formal statement that they do exist is called the Hadamard conjecture. This conjecture is over a century old and has proven itself to be one of the most difficult problems in mathematics.

A large number of Hadamard matrices are known to exist. The conjecture has been proven for all $4K \leq 664$ and there are only 13 cases that have not yet been shown for $4K \leq 2000$ \cite{dokovic_had}. Moreover, Theorem \ref{bighadthm} gives large classes of ENTIFs, found from Hadamard matrices, for all of these dimensions. Theorem \ref{bighadthm} along with the fact that Hadamard matrices are a well studied topic of research, which have yet to be classified, illustrates why classifying ENTIFs are similarly complicated. See \cite{szollosi_had} for an in-depth discussion on Hadamard matrices.

\section{Equal Norm, Tight, Integer Frames in Two and Three Dimensions}

This section addresses when ENTIFs exist in two and three dimensions. The question of existence in two dimensions is answered entirely, but only partially answered in three dimensions.

In order to obtain a full spark frame in the two dimensional case, the following result concerning the number of representations of an integer as the sum of two squares is needed.

\begin{lem} \cite[Ch.~XV]{sqrs_book} \label{sqrsthm}
Let $n = 2^{a_0}p_1^{2a_1} \cdots p_r^{2a_r} q_1^{b_1} \cdots q_s^{b_s}$, where the $p_i$'s are prime numbers of the form $4x - 1$ for $i=\{1,\cdots,r\}$, the $q_j$'s are prime numbers of the form $4x+1$ for $j=\{1,\cdots,s\}$, and ${a_i},{b_j} \in \mathbb{Z}$, for $i=\{1,\cdots,r\}$ and for $j=\{1,\cdots,s\}$. If
\[
B = (b_1 +1)(b_2+1)\cdots(b_s+1),
\]
then the number of distinct representations of $n$ as the sum of two unequal squares, ignoring order, is given by
\[
N_s(n) = \left\{\begin{array}{l@{\qquad}l}\vspace{.2cm}
\frac{B}{2} & \mbox{if $B$ is even}\\
\frac{B-1}{2} & \mbox{if $B$ is odd}
\end{array}\right.
\]
\end{lem}

As an application of this lemma, we will show that there exists a full spark, ENTIF in $\mathcal{H}_2$ with $2N$ elements for all positive integers $N >0$.

\begin{thm} \label{2N}
There exists a full spark, ENTIF in $\mathcal{H}_2$ with $2N$ elements for all positive integers $N$.
\end{thm}

\begin{proof}
Taking $n = c^2 = 5^{2N}$ (and hence $q_1=5$ and $b_1=2N$) in Lemma \ref{sqrsthm} implies that $c^2$ has $N$ distinct representations as a sum of two unequal squares, ignoring order. Hence, there exists distinct pairs $a_i, b_i \in \mathbb{Z}$ for $i=\{1,\cdots,N\}$ such that
\[
c^2 = a_1^2 + b_1^2 = \cdots = a_N^2 + b_N^2.
\]
If $A$ is the $2 \times 2N$ matrix given by
\[
A = \left[\begin{array}{rr@{\quad \cdots \quad}rr} a_1 & b_1 & a_N & b_N \\ b_1 & -a_1 & b_N & -a_N \end{array}\right],
\]
then $A$ clearly represents an ENTIF and it is full spark since each representation of $c^2$ is distinct.
\end{proof}

\begin{cor} \label{2MN}
There exists an ENTIF in $\mathcal{H}_{2M}$ with $2MN$ elements for any positive integers $M$ and $N$.
\end{cor}

\begin{proof}
Let $A$ be a $2 \times 2N$ matrix representing an ENTIF frame in $\mathcal{H}_2$ (Theorem \ref{2N} guarantees that one exists for all positive integers $N$). Let $B$ be the $2M \times 2MN$ block diagonal matrix $B = \diag(A,\ldots,A)$ obtained by adjoining $M$ copies of $A$ together as described in Proposition \ref{adjoinprop2}. Then $B$ represents a $2MN$ element ENTIF frame in $\mathcal{H}_{2M}$.
\end{proof}

\begin{rmk}
For $M>1$, the frame $B$, obtained in the proof of Corollary \ref{2MN}, is full spark only when $N = 1$, whence the frame is a basis.
\end{rmk}

We have seen that there exists a $2N$ element full spark, ENTIF in $\mathcal{H}_2$ for all positive integers $N$; however, this is not the case when the frame has $2N+1$ elements for any positive integer $N$. In fact, there does not exist any ENTIFs in $\mathcal{H}_2$ with an odd number of elements. To prove this fact, we need to carefully examine the parities of two sets of integers which square sum to the same number.

\begin{dfn}
Let $n,m \in \mathbb{N}$ and set $p=m+n$. A set of integers $(a_i)_{i=1}^p$ has \emph{parity} $[m,n]$ if $m$ integers in $(a_i)_{i=1}^p$ are even and $n$ integers in $(a_i)_{i=1}^p$ are odd.
\end{dfn}

\begin{prop}\label{oddsqrs}
Let $\{a_i\}_{i=1}^N$ and $\{b_j\}_{j=1}^M$ be integers satisfying
\[ \sum_{i=1}^N a_i^2 = \sum_{j=1}^Mb_j^2,\]
and let
\[ I = \{ 1\le i \le N: a_i \mbox{ is odd}\},\mbox{ and }J=\{1\le j \le M:b_j\mbox{ is odd}\}.\]
Then $|I|-|J|$ is divisible by 4.
\end{prop}

\begin{proof}
First note that
\[ \sum_{i\in I}a_i^2 + \sum_{i\in I^c}a_i^2 = \sum_{j\in J}b_j^2 + \sum_{j\in J^c}b_j^2.\]
and hence rearranging gives
\[ \sum_{i\in I}a_i^2-\sum_{j\in J}b_j^2 = \sum_{j\in J^c}b_j^2 - \sum_{i\in I^c}a_i^2.\]
Since all terms on the right hand side are squares of even integers, we have that $\sum_{j\in J^c}b_j^2 - \sum_{i\in I^c}a_i^2$ is divisible by 4. Next, since all terms on the
left hand side are squares of odd integers, then $\sum_{i\in I}a_i^2-\sum_{j\in J}b_j^2$ is divisible by 4 if and only if $|I|-|J|$ is divisible by 4.
\end{proof}

\begin{cor}\label{parity}
Suppose that $A = \{a_i\}_{i=1}^M$ and $B = \{b_i\}_{i=1}^M$ satisfy
\[ \sum_{i=1}^M a_i^2 = \sum_{i=1}^M b_i^2.\]
If the parity of $A$ is $[m,M - m]$, then the parity $B$ is $[m+4k,M-m-4k]$ for some integer $k$.
\end{cor}

\begin{proof}
The proof follows from Proposition \ref{oddsqrs}.
\end{proof}

Now we can show that ENTIFs with an odd number of elements do not exist in $\mathcal{H}_2$.

\begin{thm}\label{odd2dim}
An ENTIF with an odd number of elements does not exist in $\mathcal{H}_2$.
\end{thm}

\begin{proof}
Suppose by way of contradiction that there exists an ENTIF, $A$, in $\mathcal{H}_2$ with $2N+1$ elements, for some $N \in \mathbb{N}$. Note that if $A$ consisted of all even elements, then we could factor out the largest common factor of $2^k$ from each element of $A$, for some $k \in \mathbb{N}$, and we will be left with $2^k \hat{A}$, where $\hat{A}$ has at least one odd element. So without loss of generality, we may assume that $A$ has at least one odd element. Furthermore, observe that if both rows contain all odd elements, then the inner product of the rows cannot be zero, which contradicts our assumption that the rows of any integer frame must be orthogonal. Hence, $A$ must have at least one even element and at least one odd element.

Therefore, since the square sums of the columns must be equal then Corollary \ref{parity} implies that each column has parity $[1,1]$. Hence the total number of odd elements in $A$ is $2N+1$, an odd number. However, since the square sums of the rows must also be equal, then Corollary \ref{parity} also implies that if $s_1$ is the number of odd elements in the first row, then $s_1 - 4k$ is the number of odd elements in the second row for some integer $k$. Thus, the total number of odd elements in $A$ is $2(s_1 - 2k) $, which is an even number, hence a contradiction is met and such an  $A$ cannot exist.
\end{proof}

So far, we have fully classified ENTIFs in $\mathcal{H}_2$ and we would similarly like to be able to fully classify ENTIFs in $\mathcal{H}_3$. However, the three dimensional case has further complications and hence only a partial classification is obtained. First, it is shown that ENTIFs having a number of vectors that is a multiple of three or a multiple of four exist in three dimensions.

\begin{thm} \label{3dim}
For any positive integer $N$, there exists an ENTIF in $\mathcal{H}_3$ with $3N$ elements and there exists an ENTIF in $\mathcal{H}_3$ with $4N$ elements.
\end{thm}

\begin{proof}
Let $A$ be any $3\times 3$ integer matrix whose columns form an orthonormal basis for $\mathcal{H}_3$. Such matrices exist in abundance by first finding one with rational entries and then multiplying by the common denominator. However, one can simply choose $A$ to be the $3\times 3$ identity matrix. Then the matrix $[A, \cdots, A]$ obtained by adjoining $N$ copies of $A$ together, as in Proposition \ref{adjoinprop1}, is an ENTIF with $3N$ elements in $\mathcal{H}_3$. The $4N$-element case is obtained in a similar manner by adjoining $N$ copies of the $4\times 4$ Hadamard matrix, as described in Section \ref{sec_oldnew}.
\end{proof}

\begin{cor}
For any positive integers $M$ and $N$, there is an ENTIF in $\mathcal{H}_{3M}$ with $3MN$-elements and there is an ENTIF in $\mathcal{H}_{3M}$ with $4MN$-elements.
\end{cor}

\begin{proof}
Redefine the matrix $A$ in Corollary \ref{2MN} to be a  $3 \times 3N$ matrix or a $3 \times 4N$ matrix representing an ENTIF in $\mathcal{H}_3$, which is guaranteed by Theorem \ref{3dim}. Then the proof follows from the proof of Corollary \ref{2MN} where $B$ is now redefined to be a $3M \times 3MN$ block diagonal matrix, or a $3M \times 4MN$ block diagonal matrix, respectively.
\end{proof}

\begin{rmk}
Unfortunately, for any $p, q \in \mathbb{N}$, we cannot adjoin $p$ copies of a 3-element ENTIF with  $q$ copies of a 4-element ENTIF to get new ENTIFs in $\mathcal{H}_3$, because the square norms of their columns can never be the same.
\end{rmk}

Next, necessary conditions for when a matrix of size $3 \times (2N + 1)$ represents an ENTIF is given, which will lead to proving that an ENTIF with five elements in three dimensions does not exist.

\begin{thm}\label{big3dim}
If $N$ is an integer with $N \geq 2$ such that $\gcd(2N+1,3) = 1$ and $A$ is a $3 \times (2N+1)$ matrix which represents an ENTIF in $\mathcal{H}_3$, then the parity of each column must be $[2,1]$ and the number of odds in the $i^{th}$ row is of the form $4m_i + k$ with $0 \leq k < 4$. Therefore,
\begin{eqnarray}
\label{1} 4(m_1 + m_2 + m_3) + 3k = 2N+1
\end{eqnarray}
must hold. Furthermore, $4m_i+k,4m_j+k \leq N$ for some $1\leq i\neq j \leq 3$. 
\end{thm}

\begin{proof}
As in the proof of Theorem \ref{odd2dim}, it may be assumed without loss of generality that $A$ has at least one even entry and at least one odd entry.

First, consider the case in which two rows, $R_1$ and $R_2$, of $A$ have $0 \leq s_1 \leq N$ and $0 \leq s_2 \leq N$ even entries, respectively, and let $R_3$ represent the remaining row of $A$. At least one of $s_1$ or $s_2$ is nonzero since both rows having all odd elements would imply that the two rows are not orthogonal. Also, Corollary \ref{parity} implies that each column of $A$ has parity $[1,2]$ since at least one column has two odds by the assumption that $s_i \leq N$ for $i=\{1,2\}$. That is,
up to reordering the columns and/or rows, we are in the case where the frame matrix is of the form
$$
\left[
\begin{array}{cccccccc}
e & \cdots & e & o & \cdots & o & o \cdots o \\
o & \cdots & o & e & \cdots & e & o \cdots o \\
o & \cdots & o & o & \cdots & o & e \cdots e
\end{array}
\right]
$$
where $e$ symbolizes an even integer, $o$ symbolizes an odd integer and there are $s_1$ even entries in row one ($R_1$), $s_2$ even entries in row two ($R_2$) and $2N+1-s_1-s_2$ even entries in row 3 ($R_3$).

Furthermore, since the elements of $R_1$ and $R_2$ both square sum to the same number due to $A$ being a tight frame, then Corollary \ref{parity} also gives $s_2 = s_1 + 4k$ for some integer $k$. Hence, $R_3$ must have $s_1 + s_2 = 2s_1 + 4k$ odd entries and $2N- 2s_1 - 4k + 1=2\left(N-s_1-2k\right) +1$ even entries due the parity restriction of the columns. Now, since $R_3$ has an odd number of even entries and $A$ is tight, then by Corollary \ref{parity} we see that $s_1$ and $s_2$ must also be odd numbers because they possibly differ from the number of even elements in $R_3$ by a factor of four. However, taking the inner product of $R_1$ and $R_2$ gives the sum of $2(s_1 + 2k)$ even numbers and $2(N-s_1-2k)+1$ odd numbers, which must be odd. That is, the inner product cannot be zero, yielding a contradiction.

Next consider the case that two rows $R_1$ and $R_2$ have $0 \leq s_1 \leq N$ and $0 \leq s_2 \leq N$ odd entries, respectively. Then the parity of each column must be $[2,1]$ since at least one column has two even entries. Corollary \ref{parity} implies each row has $4m_i + k$ odds and equation (1) is obtained by summing the number of odds in all rows.
\end{proof}

\begin{cor}
There does not exist a five element ENTIF in $\mathcal{H}_3$.
\end{cor}

\begin{proof}
If such an ENTIF did exist, then from Theorem \ref{big3dim} there would exist integers $m \geq 0$ and $0\leq k \leq 3$ satisfying $4m + 3k = 5$. However, by substituting in $k = 0,1,2,3$, it is immediate that no such numbers exist and so a contradiction is met.
\end{proof}

Theorem \ref{big3dim} does not give a contradiction for any number of elements larger than five. For instance, there may exist an ENTIF represented by a $3 \times 7$ matrix with one odd element in each of the first two rows and five odd elements in the last row.

\begin{problem}
In $\mathcal{H}_3$, does there exist an ENTIF with $N$ elements for $N= 7,10,11...$
for the cases not covered above?  When does there exist full spark ENTIFs in
$\mathcal{H}_3$?
\end{problem}

We will see throughout this paper that it is very difficult, in general, to construct ENTIFs
with an odd number of elements except in very special cases, such as the case when the dimension
of the space is odd (and in this case, multiples of the dimension are obtained)
or for some special classes of simplexes.

\begin{problem}
Is there something fundamental about $N$ being an odd integer that presents a block to
producing ENTIFs or is it just our construction methods which are limited?
\end{problem}

The last theorem presented in this section characterizes the number of odds in each row of a matrix representing an ENTIF in $\mathcal{H}_3$, based on the parity of the columns. The proof is similar to the proof for Theorem \ref{big3dim} and so it is omitted.

\begin{thm}
Suppose $N$ is an integer with $N \geq 2$ so that $\gcd(4N+2,3) = 1$ and $A$ is a $3 \times (4N+2)$ matrix representing an ENTIF. If the parity of each column is $[2,1]$, then the number of odds in each row is of the form $4m_i + 2$ and $m = m_1 + m_2 + m_3 = N - 1$. If the parity of each column is $[1,2]$, then the number of odds in each row is $4m_i$ and $m = m_1+m_2+m_3 = 2N - 1$.
\end{thm}

\section{Equal Norm, Tight, Integer Frames with $M+1$ vectors in $M$ dimensions}

This section is dedicated to fully classifying when an ENTIF with $M+1$ vectors exists in $M$ dimensions. We show that for such a frame to exist it must be an \emph{$M$-simplex}, from which the result will follow from a previously known result.

\begin{dfn}
An \emph{$M$-simplex} is a set of $M+1$ equiangular, equal norm vectors in $M$ dimensions.
\end{dfn}

This is the generalization of a tetrahedron in three dimensions.  Recall that unit norm tight frames with $M+1$ vectors
in $M$-dimensions are all {\it unitarily equivalent} \cite{CK}.  That is, there is a unitary
operator on $\mathbb{R}^M$ which takes the $M+1$ elements of one unit norm, tight frame to the $M+1$ elements of another unit norm, tight frame.

\begin{thm}\label{oursimplex}
If $A$ is an $M \times (M+1)$ matrix representing an ENTIF, then $A$ is equiangular. Thus, the columns of $A$ form an $M$-simplex with integer coordinates.
\end{thm}

\begin{proof}
First append an additional row to $A$, which is orthogonal to and has the same norm as all rows of $A$. Call this new $\left(M+1\right) \times \left(M+1\right)$ matrix $A'$. Since the rows of $A'$ all have the same norm and the columns of $A$ all have equal norm, the added row must be of the form $[\pm a, \pm a,\ldots,\pm a,\pm a]$ for some $a \neq 0$. By possibly multiplying columns by $-1$, which does not affect the orthogonality of the rows, it may be assumed that the last row of $A'$ is $[a,a,\ldots,a,a]$.

Now, the norm squared of each row of $A$ is $(M+1)a^2$ since it must match the norm squared of the appended row. Therefore, the norm squared of each column of $A$ is $Ma^2$ because of the relationship
\begin{align*}
(M+1)c &= \sum_{j = 1}^{M+1}c = \sum_{j=1}^{M+1}\sum_{i = 1}^M A^2_{ij}\\
&= \sum_{i=1}^M \sum_{j=1}^{M+1} A^2_{ij} = \sum_{i=1}^M d = M d,
\end{align*}
where $A_{ij}$ is the entry of $A$ in the $i^{th}$ row and $j^{th}$ column, $c$ is the equal norm squared and $d$ is the tightness factor squared. Furthermore, the columns of $A'$ must be orthogonal since $A'$ is a multiple of a unitary. Therefore, the inner product of any two columns of $A$ is $-a^2$ and so $A$ is equiangular.
\end{proof}

The full classification for when an $M$-simplex with integer coordinates exists was first
proved by I.J. Schoenberg in \cite{regsimporig} and was stated in a clearer fashion by I.G. Macdonald  in \cite{regsimpian} as follows.

\begin{thm}\label{simplexthm} \cite{regsimpian}
There exists a regular $M$-simplex in $\mathbb{R}^M$ with vertices in $\mathbb{Z}^M$ if and only if $M+1$ is the sum of $1$, $2$, $4$ or $8$ odd squares.
\end{thm}

\begin{rmk}
Theorem \ref{simplexthm} along with Theorem \ref{oursimplex} imply that an $M+1$ element ENTIF in $M$ dimensions does not exist for
\[
M = 2, 4, 5, 10, 12, 13, 14, 16, 18, 20, 21, 22, 26,\ldots.
\]
\end{rmk}

Next, an explicit construction of an ENTIF for the allowable values of $M$ is given. Note that it is equivalent to constructing a regular $M$-simplex with vertices in $\mathbb{Q}^M$. The ideas presented are mostly due to R. Chapman \cite{rchap}.

Define $m = M+1$ and let $e_1,\ldots,e_m$ be the standard orthonormal basis of $\mathbb{Q}^m$. Put $v = e_1 + \cdots + e_m$. The main idea of the construction is to find a linear operator $S$ on $\mathbb{Q}^m$ so that $S = T/\sqrt{m}$ and satisfies $Sv = e_m$, where $T$ is an orthogonal matrix. Such an $S$ preserves inner products and furthermore the set $\{S e_j\}_{j=1}^m$ forms another orthogonal set in which the $m$-th coordinate of $S e_j$ is $1/m$ for all $1 \leq j \leq m$. Therefore, removing the last row of the matrix representation with respect to the standard orthonormal basis of $S$ gives an $M+1$ element ENTIF in $M$ dimensions.

To construct such an $S$, it is enough to find a linear operator $U:\mathbb{Q}^m \to \mathbb{Q}^m$ so that $U = Q/\sqrt{m}$ for some orthogonal operator $Q$ and then compose $U$ with the reflection $R$, the hyperplane with normal vector $$\dfrac{Uv - e_m}{\| U v - e_m \|}.$$
That is, $S = R\circ U$ and so $Sv = R(Uv) = e_m$ as required.

In the case that $m$ is a perfect square, define $U x = x/\sqrt{m}$. If $m$ is the sum of $k = 2,4$, or $8$ odd squares, such as $m=a^2 + \cdots +h^2$, then let $Ux = A_k x/m$ where $A_k$ is the block diagonal matrix having $E_k$ down the diagonal $m/k$ times and where the $E_k$ are defined as
\[
E_2 = \left[\begin{array}{rr} a & -b \\ b & a \end{array}\right] ,
\]
\[
E_4 = \left[\begin{array}{rrrr} a & b & c & d \\ -b & a & -d & c \\ -c & d & a & -b \\ -d &-c & b & a \end{array}\right] ,
\]
\[
E_8 = \left[\begin{array}{rrrrrrrr}
a & b & c & d & e & f & g & h \\
-b & a & -d & c & -f & e & -h & g\\
e & -f & g & -h & -a & b & -c & d\\
-f & -e & h & g & b & a & -d & -c\\
-d & -c & b & a & -h & -g & f & e\\
c & -d & -a & b & -g & h & e & -f\\
g & -h & -e & f & c & -d & -a & b\\
-h & -g & -f & -e & d & c & b & a
\end{array}\right].
\]
The operators given in each case are easily checked to have the described properties.

In the construction above, an $(M+1)\times (M+1)$ rational unitary matrix having a row with all entries being the same modulus was constructed. Any such matrix yields an ENTIF with $M+1$ elements in $M$ dimensions by removing the constant modulus row.  An identical proof technique as in the proof of Theorem \ref{oursimplex} combined with Theorem \ref{simplexthm} immediately implies these types of matrices exist if and only if $M+1$ is the sum of $1$, $2$, $4$, or $8$ odd squares, which is stated in the following Theorem.

 \begin{thm}
 There is an equal norm tight integer frame with $(M+1)$-elements in $\mathcal{H}_M$ if and only if $M+1$ is
 the sum of $1,2,4,$ or $8$ odd squares.
 \end{thm}

\section{General Equal Norm, Tight, Integer Frames}

This section includes all of the remaining results concerning ENTIFs in a general dimension. The main result in this section gives a way to adjoin two ENTIFs to obtain an ENTIF with $N$ elements for all large enough $N$. In order to obtain this result, a basic number theory result is needed.

\begin{lem}\cite{integers} \label{integers}
If $a, b \in \mathbb{N}$ such that gcd$\left(a, b\right) = 1$, then for all integers $m \geq \left(a - 1\right)\left(b - 1\right)$, there is exactly one pair of nonnegative integers $p$ and $q$ such that $q < a$ and $m = pa + qb$.
\end{lem}

\begin{cor} \label{adjoinnums}
If $a,b \in \mathbb{Z}$ and $g$ is defined to be $g: = \mbox{gcd}(a,b)$, then for every integer $m \geq (a/g - 1)(b/g - 1)$ there exist nonnegative integers $p$ and $q$ so that $gm = pa + qb$.
\end{cor}

\begin{proof}
Note that $\gcd(a/g,b/g) = 1$, so Lemma \ref{integers} applies. Hence, there exist nonnegative integers $p$ and $q$ such that $m = p(a/g) + q(b/g)$.
\end{proof}

Combining Lemma \ref{integers} and Corollary \ref{adjoinnums} yields a fundamental result which states that if we can construct two
ENTIFs in $\mathcal{H}_M$ such that the number of vectors in the two frames are relatively prime with the same equal norm constant,
then we can construct ENTIFs with $N$-elements for all large $N$.

\begin{thm}\label{gcdframethm}
Suppose $A$ and $B$ represent ENTIFs in $\mathcal{H}_M$ with $N_1$ and $N_2$ elements, respectively, such that $A$ and $B$  have the same equal norm constant. If $K = \mbox{gcd}(N_1,N_2)$, then there is a $KN$ element ENTIF in $\mathcal{H}_M$ for all $N \geq (N_1/K-1)(N_2/K-1)$.
\end{thm}

\begin{proof}
If $N \geq (N_1/K - 1)(N_2/K - 1)$, then Corollary \ref{adjoinnums} implies the existence of nonnegative $c_N$ and $d_N$ such that that $KN = c_N \cdot N_1 + d_N \cdot N_2$. Therefore, Proposition \ref{adjoinprop1} implies that the block matrix
$$
[A, \ldots, A, B, \ldots, B],
$$
where $A$ appears $c_N$ times and $B$ appears $d_N$ times is an ENTIF in $\mathcal{H}_M$ with $KN$ elements.
\end{proof}

Theorem \ref{gcdframethm} leads to a number of corollaries implying the existence of ENTIFs.

\begin{cor}\label{oddsinM2dim}
If $M\geq 3$ is an odd integer and $K$ is the smallest integer such that $2^K \geq M^2$, then there is an ENTIF with $N$ elements in $\mathcal{H}_{M^2}$ for all $N  \geq (M^2-1)(2^K-1)$.
\end{cor}

\begin{proof}
The matrix $A = M \cdot I_{M^2 \times M^2}$ is an ENTIF with vectors having
square norms $M^2$. Furthermore, an $M^2 \times 2^K$ frame matrix $B$ which represents an ENTIF may be obtained from a $2^K \times 2^K$ Hadamard matrix (see Section \ref{sec_oldnew}) where the square norms of the columns of $B$ are also $M^2$. Since $\mbox{gcd}(M^2,2^K) = 1$, Theorem \ref{gcdframethm} gives the desired result.
\end{proof}

\begin{cor}
If $P$ is an odd integer and $M = 2P$, and $K\geq 2$ is the smallest integer such that $2^K \geq M^2$, then there is an ENTIF with $4N$ elements in $\mathcal{H}_{M^2}$ dimensions for all $N \geq (P^2 - 1)(2^{K-2} - 1)$.
\end{cor}

\begin{proof}
Choose $A$ and $B$ in exactly the same way as in the proof of Corollary \ref{oddsinM2dim}. Since $\gcd(M^2,2^{K}) = 4$, Theorem \ref{gcdframethm} gives the result.
\end{proof}

The next corollary is particularly interesting because it eliminates the necessity of knowing each Hadamard matrix before being able to construct certain ENTIF. It proves that if we have knowledge of two consecutive Hadamard matrices then we know a large class of ENTIFs exist.

\begin{cor}\label{adjoinhadamards}
If both $4N \times 4N$ and $4(N+1)\times 4(N+1)$ Hadamard matrices exist for $4N \geq M$, then for all $K \geq N(N-1)$ there is a $4K$ element ENTIF in $\mathcal{H}_M$.
\end{cor}

\begin{proof}
Since $\gcd(4N,4(N+1)) = 4$, Theorem \ref{gcdframethm} implies a $4K$ element ENTIF in $\mathcal{H}_M$ exists for $K \geq (4N/4 - 1)(4(N+1)/4-1) = N(N-1)$.
\end{proof}

The next example demonstrates the usefulness of Corollary \ref{adjoinhadamards}.

\begin{example}
Since $8\times 8$ and $12 \times 12$ Hadamard matrices exist, there are $4K$ element ENTIFs in $M \leq 8$ dimensions for all $K \geq 2$. Since only 13 Hadamard matrices are left to be shown to exist for all $4N \leq 2000$ (see Section \ref{sec_oldnew}), Corollary \ref{adjoinhadamards} gives a vast amount of ENTIFs in a large number of dimensions.
\end{example}

Next, we prove that there exists an ENTIF in $\mathcal{H}_5$ with an even number of elements for almost every positive even integer.

\begin{cor}\label{evensin5dim}
For every $N \geq 12$, there is a $2N$ element ENTIF in $\mathcal{H}_5$.
\end{cor}

\begin{proof}
Let $a$\ be any nonzero integer and let $b = 2a$. Then the $5 \times 8$ matrix
\[
A = \left[\begin{array}{rrrrrrrr}
a & a & a & a & a & a & a & a\\
b & -b & 0 & 0 & 0 & 0 & 0 & 0\\
0 & 0 & b & -b & 0 & 0 & 0 & 0\\
0 & 0 & 0 & 0 & b & -b & 0 & 0\\
0 & 0 & 0 & 0 & 0 & 0 & b & -b
\end{array}\right]
\]
and the $5 \times 10$ matrix
\[
B = \left[\begin{array}{rrrrrrrrrr}
a & -b & 0 & 0 & 0 & 0 & 0 & 0 & a & -b\\
b & a & a & -b & 0 & 0 & 0 & 0 & 0 & 0\\
0 & 0 & b & a & a & -b & 0 & 0 & 0 & 0\\
0 & 0 & 0 & 0 & b & a & a & -b & 0 & 0\\
0 & 0 & 0 & 0 & 0 & 0 & b & a & b & a
\end{array}\right]
\]
represent ENTIFs having the same equal norm squared, $a^2 + b^2$, so that Theorem \ref{gcdframethm} gives a $2N$-element ENTIF in $\mathcal{H}_5$ for all $N \geq 12$.
\end{proof}

\begin{rmk}
By Theorem \ref{simplexthm}, a six element ENTIF does not exist in five dimensions. Due to  Section \ref{sec_oldnew} and Theorem \ref{bighadthm}, since the $8 \times 8$, $12\times 12$, $16 \times 16$ and $20 \times 20$ Hadamard matrices exists, then there exist ENTIFs in $\mathcal{H}_5$ with 8, 12, 16, and 20 elements. Also, adjoining two copies of the $5 \times 5$ identity matrix, as in Proposition \ref{adjoinprop1}, yields a 10-element ENTIF in $\mathcal{H}_5$. Lastly, Corollary \ref{evensin5dim} proves the existence of ENTIFs with an even number of vectors in $\mathcal{H}_5$ for all even integers $N \geq 24$. Therefore, the only even element ENTIFs in $\mathcal{H}_5$ for which the existence is unknown are those with $N = 14, 18$ and $22$ elements.
\end{rmk}

The last theorem of this section gives the existence of $4N^2$ and $8N^2$ element ENTIFs, from which Theorem \ref{gcdframethm} can be applied to obain even more.

\begin{thm}\label{gensqrthm}
If $N$ is a positive integer, then there exists an ENTIF with
\begin{enumerate}
\item $4N^2$ vectors in $N^2 + 1$ dimensions
\item $4N^2$ vectors in $2N^2 + 1$ dimensions
\item $4N^2$ vectors in $3N^2 + 1$ dimensions
\item $8N^2$ vectors in $4N^2 + 1$ dimensions
\item $8N^2$ vectors in $4N^2 + 2$ dimensions
\end{enumerate}
\end{thm}

\begin{proof} For (1)-(3), let $b$ be a nonzero integer and $a = Nb$.

(1) For each $1 \leq j \leq N^2$, let $B_j$ be the $(N^2+1) \times 4$ matrix $[b,b,b,-b]$ as its first row, $[a, a, -a, a]$ as its $(j+1)$ row, and all other rows having zero entries. If $A$ is the $(N^2+1)\times 4N^2$ matrix given by $A = \left[B_1,\ldots,B_{N^2}\right]$, then the choice of $a$ and $b$ ensure that $A$ is the desired ENTIF. Note that the equal norm squared is $(N^2+1)b^2$.

(2)  For each $1 \leq j \leq N^2$, let $B_j$ be the $(2N^2+1)\times 4$ matrix having $[b,b,b,-b]$ as its first row, $[a,a,-a,a]$ as its $2j$ row, $[a,-a,a,a]$ as its $(2j+1)$ row, and all other rows having zero entries. If $A$ is the $(2N^2+1)\times 4N^2$ matrix given by $A = \left[B_1,\ldots,B_{N^2}\right]$, then the choice of $a$ and $b$ ensures that $A$ is the desired ENTIF. Note that the equal norm squared is $(2N^2+1)b^2$.

(3) For each $1 \leq j \leq N^2$, let $B_j$ be the $(3n^2+1)\times 4$ matrix having $[b,b,b,-b]$ as its first row, $[a,a,-a,a]$ as its $(3j-1)$ row, $[a,-a,a,a]$ as its $3j$ row, $[-a,a,a,a]$ as its $(3j+1)$ row, and all other rows having zero entries. If $A$ is the $(3N^2+1)\times 4N^2$ matrix given by $A = \left[B_1,\ldots,B_{N^2}\right]$, then the choice of $a$ and $b$ ensures that $A$ is the desired ENTIF. Note that the equal norm squared is $(3N^2+1)b^2$.

For (4) and (5), let $b$ be a nonzero integer and $a = 2Nb$.

(4) For each $1 \leq j \leq N^2$, let $B_j$ be the $(4N^2+1)\times 8$ matrix having $[b,b,b,b,b,b,b,b]$ as its first row, having $[a,-a,0,0,0,0,0,0]$ as its $4j - 2$ row, $[0,0,a,-a,0,0,0,0]$ as its $4j-1$ row, $[0,0,0,0,a,-a,0,0]$ as its $4j$ row, $[0,0,0,0,0,0,a,-a]$ as its $4j+1$ row, and zero entries in all other rows. If $A = \left[B_1,\ldots,B_{N^2}\right]$, then the choice of $a$ and $b$ ensures that $A$ is the desired ENTIF. Note that the equal norm squared is $(4N^2+1)b^2$.

(5) For each $1 \leq j \leq 2 N^2$, let $B_j$ be the $(4N^2+2)\times 4$ matrix having $[b,b,b,b]$ as its first row, $[b,-b,b,-b]$ as its second row, having $[a,0,-a,0]$ as its $2j + 1$ row, $[0,a,0,-a]$ as its $2j+2$ row, and zero entries in all other rows. If $A = \left[B_1,\ldots,B_{2N^2}\right]$, then the choice of $a$ and $b$ ensures that $A$ is the desired ENTIF. Note that the equal norm squared is $2b^2(1+2N^2)$.
\end{proof}

\begin{example}
Theorem \ref{gensqrthm} says that there exists ENTIFs with:
\begin{enumerate}
\item 4 vectors in $\mathcal{H}_2$, 16 in $\mathcal{H}_5$, 36 in $\mathcal{H}_{10}$, 64 in $\mathcal{H}_{17}$, \ldots
\item 4 vectors in $\mathcal{H}_3$, 16 in $\mathcal{H}_{9}$, 36 in $\mathcal{H}_{19}$, 64 in $\mathcal{H}_{33}$, \ldots
\item 4 vectors in $\mathcal{H}_4$, 16 in $\mathcal{H}_{13}$, 36 in $\mathcal{H}_{28}$, 64 in $\mathcal{H}_{49}$, \ldots
\item 8 vectors in $\mathcal{H}_5$, 32 in $\mathcal{H}_{17}$, 72 in $\mathcal{H}_{37}$, 128 in $\mathcal{H}_{65}$, \ldots
\item 8 vectors in $\mathcal{H}_6$, 32 in $\mathcal{H}_{18}$, 72 in $\mathcal{H}_{38}$, 128 in $\mathcal{H}_{66}$, \ldots
\end{enumerate}
Furthermore, these ENTIFs can be adjoined to obtain multiplies of the given number of vectors.
\end{example}

\begin{rmk}
One can construct an $(N^2 + 1) \times 2(N^2+1)$ matrix in a similar fashion as matrix $B$ in the proof of Corollary \ref{evensin5dim} and adjoin it with the matrix in Theorem \ref{gensqrthm}(1) to obain an ENTIF in $N^2 + 1$ dimensions with $4K$ elements for all $K \geq N^2(N^2-1)$ by Theorem \ref{gcdframethm}. One can also do the same in $4N^2+1$ dimensions to obtain a $2K$ element ENTIF with $2K$ elements for all $K \geq 4N^2(4N^2-1)$.
\end{rmk}

\section{Removing Either the Equal Norm Or Tightness Assumption}

Only ENTIFs have been considered so far.
As we have seen, these can be quite difficult to construct.  So in this section,
we address
 the question of what can be obtained if one of the assumptions that the frame is equal norm or tight  is removed. In either case, it will be shown that an integer frame of any size in any dimension may be obtained.

\begin{thm} \label{en}
If $M$ and $N$ are positive integers satisfying $N \geq M$, then there is an equal norm integer frame with $N$ elements in $\mathcal{H}_M$.
\end{thm}

\begin{proof}
For all $1\leq i \leq M$, let $A_{i}$ be the $M \times i$ matrix formed by the first $i$ columns of the identity matrix $I_{M\times M}$. Write $N = cM + k$ for some integers $c \geq 0$ and $0 \leq k < M$. If $k = 0$, then the block matrix $C = [A_M \cdots A_M]$ where $A_M$ is repeated $c$ times is an equal norm (tight) integer frame with $N$ elements. If $k > 0$, then the block matrix $C = [A_M \cdots A_M \mbox{ } A_k]$ where $A_M$ is repeated $c$ times is a desired equal norm integer frame.
\end{proof}

Before proving that tight integer frames exist with any number of elements in any dimension,
the following number theoretic result is needed.

\begin{lem}\label{arbnumsqrs}
For every positive integer $k$, there exists a nonzero integer $s$ such that $s^2$ can be written as a sum of $i$ nonzero squares for all $1 \leq i \leq k$.
\end{lem}

\begin{proof}
First recall the well-known Euclid's formula, which states that if $m$ and $n$ are positive integers with $m > n$, then
\[
a = m^2 - n^2, \quad b = 2 m n, \quad c = m^2 + n^2
\]
forms a Pythagorean triple, i.e., $a^2 +b ^2 = c^2$. Suppose that $m_0$ and $n_0$ are odd integers with $m_0 > n_0$ and let $(a_0,b_0,c_0)$ be the Pythagorean triple formed by $m_0$ and $n_0$ as given by Euclid's Formula. Since $m_0$ and $n_0$ are both odd, $c_0 = 2 \cdot m_1$ for some odd integer $m_1$. Letting $n_1 = 1$ gives another Pythagorean triple $(a_1,b_1,c_1)$ generated by $m_1$ and $n_1$ in which $c_1^2=a_1^2+b_1^2$ and $b_1 = 2 m_1 \cdot  n_1 = c_0$. Thus
\[
c_1 ^2 = a_1^2 + b_1^2 = a_1^2 + c_0^2 = a_1^2 + a_0^2 + b_0^2.
\]
This process may be continued to find a number $c_{i-2}$, such that $c_{i-2}^2$ is the sum of $3 \leq i \leq k$ squares. This follows because in each step $c_{i-3}$ is always of the form $2 m_{i-2}$ for some odd integer $m_{i-2}$ and so $b_{i-2} = c_{i-3}$ with $n_{i-2} = 1$.
\end{proof}

\begin{thm} \label{tight}
If $M$ and $N$ are positive integers satisfying $N \geq M$, then there is a tight integer frame with $N$ elements in $\mathcal{H}_M$.
\end{thm}

\begin{proof}
If $M$ is even, let $k = (M-2)/2$. Let $p$ be a nonzero integer such that $p^2$ can be written as a sum of $i$ nonzero squares for all $1 \leq i \leq N-2k-1$, which exists by Lemma \ref{arbnumsqrs}. Write
\[
p^2 = a^2 + b^2 = a_1^2 + \cdots + a_{N-2k-1}^2
\]
for some $a,b,a_i \in \mathbb{Z}$ and define
\[
A_i = \left[ \begin{array}{rr} a & b \\ b & -a \end{array}\right]
\]
for all $1 \leq i \leq k$. Define $A$ to be the $M \times N$ matrix given by
\[
A = \left[ \begin{array}{ccccccc}
A_1 & & & & \\
 & \ddots & & &\\
 & & A_k & & \\
 & & & a_1 & \cdots & a_{N-2k-1} & \\
 & & & & & & p
\end{array}\right]
\]
 where all of the empty entries in $A$ are $0$. Then $A$ is a tight integer frame with $N$ elements in $M$ dimensions with tightness factor $p$.

If $M$ is odd, let $k = (M-1)/2$ and let $p$ be a nonzero integer such that $p^2$ can be written as a sum of $i$ nonzero squares for all $1 \leq i \leq N-2k$. Write
\[
p^2 = a^2 + b^2 = a_1^2 + \cdots + a_{N-2k}^2
\]
for some $a,b,a_i \in \mathbb{Z}$ and define $A_i$ as above for all $1 \leq i \leq k$. Define $A$ to be the $M \times N$ matrix given by
\[
A = \left[ \begin{array}{cccccc}
A_1 & & & & \\
 & \ddots & & &\\
 & & A_k & & \\
 & & & a_1 & \cdots & a_{N-2k} \\
\end{array}\right]
\]
where all of the empty entries in $A$ are $0$. Then $A$ is a tight integer frame with $N$ elements in $M$ dimensions with tightness factor $p$.
\end{proof}

Throughout this paper, we have made numerous theorems and classifications for when ENTIFs, equal norm integer frames, and tight integer frames exist. Although we have proven that ENTIFs with an odd number of elements do not exist in $\mathcal{H}_2$, and we have proven other general statements about the existence of ENTIFs, there has yet to be a complete classification for when ENTIFs with $N$ elements exist in $\mathcal{H}_M$. We have, however, seen that there exist $N$ element tight frames and $N$ element equal norm frames in $\mathcal{H}_M$ for all $N \geq M$. This next result looks at frames which are almost ENTIFs and implies that there exists an equal norm integer frame in $\mathcal{H}_M$ with $N$ elements for $N \geq M$ that is arbitrarily close to being tight.

First, the formal definition of a frame being arbitrarily close to tight is given and then the result is stated and proved.

\begin{dfn}
A frame $\{f_i\}_{i=1}^N$ is said to be \emph{$(\varepsilon, A)$-tight} if there are constants $0 < \varepsilon < 1$ and $A > 0$ such that the lower and upper frame bounds are $(1-\varepsilon)A$ and $(1+\varepsilon)A$, respectively.
\end{dfn}

\begin{thm}\label{almosttight}
Let $M$ and $N$ be positive integers such that $N \geq M$. For any $\varepsilon > 0$ and any orthonormal basis $\beta = \{e_i\}_{i=1}^M$ for $\mathcal{H}_M$, there exists a full spark, equal norm, integer frame $F = \{f_i\}_{i=1}^N$ with respect to $\beta$ for which  $F$ is $(\varepsilon , N/M)$-tight.
\end{thm}

\begin{proof}
It is enough to show the existence of such a frame with rational coordinates. Begin by first picking a unit norm tight frame $\Psi = \{\psi_i\}_{i=1}^N$ with tight frame bound $N/M$
\cite{petesbook}. Note that $\Psi$ may not have rational coordinates.

Let $0 < \varepsilon < 1$ be given and momentarily fix a $0 < \delta < 1$, which will be chosen later. Since vectors with rational coordinates are dense in $S^{M-1}$, the unit sphere in $\mathbb{R}^M$, vectors $F_1 = \{f_i\}_{i=1}^M$ with rational coordinates can be chosen to be linearly independent and satisfy $$\| f_i - \psi_i \| \leq \frac{\delta}{\sqrt{M}}$$ for all $1 \leq i \leq M$.

Now let $\mathbb{H}_1$ be the collection of all hyperplanes in $\mathbb{R}^M$ generated by
sets of $M-1$ vectors chosen
from
$F_1$. If
\begin{align*}
C_1 = \left(\bigcup_{H \in \mathbb{H}_1} H\right)^c,
\end{align*}
then $C_1$ is also dense in $S^{M-1}$ and so
we can choose $f_{M+1}$ in $S^{M-1} \cap C_1$ with rational coordinates so that $$\| f_{M+1} - \psi_{M+1} \| \leq \frac{\delta}{\sqrt{M}}.$$ Notice that by construction the set $F_2 = \{f_i\}_{i=1}^{M+1}$ is full spark. Now choose the set of all hyperplanes $\mathbb{H}_2$ generated by $F_2$ and continue the same process until a frame $F = \{f_i\}_{i=1}^N$, having all rational coordinate vectors, is obtained.

Now to prove that $F$ is $(\epsilon,N/M)$-tight. Minkowski's inequality and the Cauchy-Schwarz inequality gives for any $x \in \mathbb{R}^M$,
\begin{align*}
\left(\sum_{i=1}^N \left| \langle x, f_i \rangle \right|^2\right)^{\frac{1}{2}} &\leq \left(\sum_{i=1}^N\left| \langle x, \psi_i \rangle \right|^2 \right)^\frac{1}{2} + \left(\sum_{i =1}^N \left| \langle x, f_i - \psi_i \rangle \right|^2 \right)^{\frac{1}{2}}\\
&\leq \sqrt{\frac{N}{M}} \|x \| + \left(\sum_{i=1}^N \|x\|^2 \|f_i - \psi_i\|^2  \right)^{\frac{1}{2}}\\
& \leq \|x\| \left[\sqrt{\frac{N}{M}} + \left(\sum_{i=1}^N \frac{\delta^2}{M}\right)^{\frac{1}{2}}\right]\\
&= \|x\|(1+\delta) \sqrt{\frac{N}{M}}
\end{align*}
proving that an upper frame bound of $F$ is $(1+\delta)^2 N/M$. Similarly, a lower frame bound of $F$ is $(1-\delta)^2 N/M$. Now choose $\delta$ so that
\[
(1 - \varepsilon) \cdot \frac{N}{M} \leq (1 -  \delta)^2 \cdot \frac{N}{M} \leq (1 + \delta)^2 \cdot \frac{N}{M} \leq (1 + \varepsilon) \cdot \frac{N}{M},
\]
showing that $F$ is $(\varepsilon,N/M)$-tight.
\end{proof}

\begin{rmk}
The frame $F$ constructed in Theorem \ref{almosttight} is not necessarily represented against the eigenbasis of its frame operator as in all previous cases in the present paper.
\end{rmk}

The proof of Theorem \ref{almosttight} relies heavily on the fact that the set of all rational coordinate points are dense in $S^{M-1}$. Unfortunately, the higher the dimension and the closer the frame is to being tight forces the need to choose numbers in which the denominators are possibly massive. That is, using the proof technique above might lead to computationally inconvenient integer frames after clearing out the denominators. See \cite{schmutz1} for more details concerning rational coordinate points on the sphere.

It is also worth noting that the technique used to prove Theorem \ref{almosttight}, is a standard argument which shows that full spark, equal norm frames are {\it dense} in the space of all equal norm frames \cite[Ch.~4]{petesbook}. It is an open problem whether
the full spark, equal norm, Parseval frames are dense in the space of all equal norm
Parseval frames.

This paper is the beginning of the study on ENTIFs. We have made numerous characterizations of integer frames throughout this paper; however, there are still a lot of interesting and important open problems, as we have seen.  For one, by adjoining matrices to get larger ones, we give up full spark.  Also, this gives frames with many repeated frame vectors, which we do not usually want in practice since this repetition gives no new information.  We believe that many of the open problems here will require a deep knowledge of number theory for their resolution.




\end{document}